\documentclass{article}

\usepackage{arxiv}
\usepackage{lipsum}
\usepackage{amsmath}
\usepackage{xr-hyper}
\usepackage{hyperref}
\usepackage{graphics,graphicx}
\usepackage{epstopdf}
\usepackage{amsfonts}
\usepackage{graphicx}
\usepackage{epstopdf}
\usepackage{algorithm}
\usepackage{algorithmic}
\usepackage[algo2e,linesnumbered,lined, ruled]{algorithm2e}
\usepackage[capitalize,nameinlink]{cleveref}
\usepackage{amssymb}
\usepackage{nicefrac}
\usepackage{bm} 
\usepackage{mathtools}
\usepackage{pgfplotstable} 
\usepackage{booktabs}
\usepackage{multicol}
\usepackage{colortbl}
\usepackage{makecell}
\usepackage{siunitx}
\usepackage{ifpdf}

\usepackage[caption=false]{subfig} 

\usepackage{amsopn}
\usepackage{amsthm}

\Crefname{figure}{Figure}{Figures}

\crefformat{equation}{\textup{#2(#1)#3}}
\crefrangeformat{equation}{\textup{#3(#1)#4--#5(#2)#6}}
\crefmultiformat{equation}{\textup{#2(#1)#3}}{ and \textup{#2(#1)#3}}
{, \textup{#2(#1)#3}}{, and \textup{#2(#1)#3}}
\crefrangemultiformat{equation}{\textup{#3(#1)#4--#5(#2)#6}}%
{ and \textup{#3(#1)#4--#5(#2)#6}}{, \textup{#3(#1)#4--#5(#2)#6}}{, and \textup{#3(#1)#4--#5(#2)#6}}

\Crefformat{equation}{#2Equation~\textup{(#1)}#3}
\Crefrangeformat{equation}{Equations~\textup{#3(#1)#4--#5(#2)#6}}
\Crefmultiformat{equation}{Equations~\textup{#2(#1)#3}}{ and \textup{#2(#1)#3}}
{, \textup{#2(#1)#3}}{, and \textup{#2(#1)#3}}
\Crefrangemultiformat{equation}{Equations~\textup{#3(#1)#4--#5(#2)#6}}%
{ and \textup{#3(#1)#4--#5(#2)#6}}{, \textup{#3(#1)#4--#5(#2)#6}}{, and \textup{#3(#1)#4--#5(#2)#6}}

\theoremstyle{plain}

\theoremstyle{definition}
\newtheorem{definition}{Definition}[section]

\ifpdf
  \DeclareGraphicsExtensions{.eps,.pdf,.png,.jpg}
\else
  \DeclareGraphicsExtensions{.eps}
\fi

\DeclareSIUnit\Dalton{Da}

\usepackage{enumitem}
\setlist[enumerate]{leftmargin=.5in}
\setlist[itemize]{leftmargin=.5in}

\pgfplotsset{compat = 1.16}

\definecolor{myGreen}{RGB}{25,142,33}

\newcommand{\R}{\mathbb{R}}

\DeclareMathOperator{\id}{id}
\DeclareMathOperator{\med}{med}
\DeclareMathOperator{\diag}{diag}
\DeclareMathOperator{\dist}{dist}
\DeclareMathOperator{\weight}{w}

\newcommand{\email}[1]{\protect\href{mailto:#1}{#1}}
\newcommand{\myvec}[1]{\bm{\mathit{#1}}}

\crefname{hypothesis}{Hypothesis}{Hypotheses}

\Crefname{ALC@unique}{Line}{Lines} 

\crefname{algocf}{alg.}{algs.}
\Crefname{algocf}{Algorithm}{Algorithms}

\title{Regularized Orthogonal Nonnegative Matrix Factorization and \texorpdfstring{$K$}{K}-means Clustering}

\author{\href{https://orcid.org/0000-0003-3424-8031}{\includegraphics[scale=0.06]{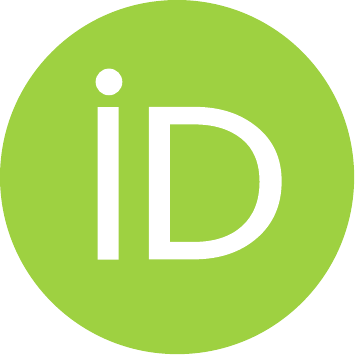}\hspace{1mm}Pascal~Fernsel}\\
        Center for Industrial Mathematics\\
        University of Bremen\\
        Bremen, Germany \\
        \texttt{\email{p.fernsel@uni-bremen.de}} \\
    \And
	\href{https://orcid.org/0000-0003-1448-8345}{\includegraphics[scale=0.06]{orcid.pdf}\hspace{1mm}Peter~Maa\ss} \\
	Center for Industrial Mathematics\\
	University of Bremen\\
	Bremen, Germany \\
	\texttt{\email{pmaass@uni-bremen.de}} \\
}

\renewcommand{\shorttitle}{Regularized ONMF and \texorpdfstring{$K$}{K}-means Clustering}

\ifpdf
\hypersetup{
  pdftitle={Regularized Orthogonal Nonnegative Matrix Factorization and K-means CLustering},
  pdfauthor={P. Fernsel},
  pdfkeywords={Orthogonal nonnegative matrix factorization, Clustering, K-means, Regularization}
}
\fi

\begin{document}

    \maketitle
    
    \begin{abstract}
        In this work, we focus on connections between $K$-means clustering approaches and Orthogonal Nonnegative Matrix Factorization (ONMF) methods. We present a novel framework to extract the distance measure and the centroids of the $K$-means method based on first order conditions of the considered ONMF objective function, which exploits the classical alternating minimization schemes of Nonnegative Matrix Factorization (NMF) algorithms. While this technique is characterized by a simple derivation procedure, it can also be applied to non-standard regularized ONMF models. Using this framework, we consider in this work ONMF models with $\ell_1$ and standard $\ell_2$ discrepancy terms with an additional elastic net regularization on both factorization matrices and derive the corresponding distance measures and centroids of the generalized $K$-means clustering model. Furthermore, we give an intuitive view of the obtained results, examine special cases and compare them to the findings described in the literature.
    \end{abstract}

    \keywords{Orthogonal nonnegative matrix factorization \and \texorpdfstring{$K$}{K}-means Clustering \and Distance function \and Centroid \and Elastic net regularization}

    \AMSClass{15A23 \and 65F22 \and 62H30}

    \section{Introduction} \label{sec:Introduction}
        Cluster analysis has been studied extensively for more than five decades by the machine learning community and has numerous applications in computer science, social science, biology, medicine and many other fields. Clustering in its classical form is an unsupervised learning technique and is a main task of exploratory data analysis. The most well-known clustering algorithms are based on either partitional or hierarchical clustering techniques. Partitional clustering methods are iterative approaches and try to optimize an objective function, while hierarchical clustering algorithms typically develop a binary tree-based data structure to obtain the desired clustering.
        
        In this work, we focus on the $K$-means clustering, which is the most widely used partitional clustering algorithm, and analyze its relationships to different \textit{Nonnegative Matrix Factorization} (NMF) models, which are a specific matrix factorization method with additional nonnegativity constraints. The clustering capability of NMF approaches are well-known throughout the literature and the study of the connections between both approaches leads to several different NMF models for clustering, which typically have advantages compared with usual $K$-means approaches. However, the corresponding works only consider non-regularized NMF models and mainly establish these connections by directly comparing the objective functions and the constraints of the optimization problems of both the NMF problem as well as the K-means model.

        This work rather focuses on obtaining relationships between regularized Orthogonal NMF (ONMF) and generalized $K$-means models by deriving the corresponding distance measures and centroids of the K-means approach via the typical alternating minimization scheme of NMF algorithms and first order conditions of the ONMF objective function. The used approaches to obtain the desired connections between both approaches are significantly easier compared to the proof techniques throughout the literature and generalize to regularized ONMF models. In this way, from a viewpoint of the field of inverse problems, the effect of the considered elastic net regularization on the obtained distance function and centroid of the K-means approach can be obtained. Furthermore, classical K-means approaches using usual distance functions typically suffer from the uniformization effect and tend to produce clusters with relatively balanced cluster sizes. However, the used technique in this work can lead to generalized K-means models which do not suffer from the uniformization effect.

        This paper is organized as follows. \cref{sec:Background and Related Work} is devoted to the basics of $K$-means clustering, ONMF as well as their relationships and gives an overview of the related works. \cref{sec:Regularized ONMF Model and Generalized K-means} introduces the considered ONMF models of this work and discusses the relationships to the corresponding generalized $K$-means models. After a short description of some basic properties of the considered framework in \cref{subsec:Alternating Minimization and Separability}, the following Sections \ref{subsec:Classical Discrepancy} and \ref{subsec:l1 Discrepancy} contain a detailed derivation of the distance functions and centroids of the generalized $K$-means models based on the considered ONMF problem for different discrepancy terms. Furthermore, \cref{subsec:Evaluation and Special Cases} shortly discusses some special cases and compare the obtained findings to known results in the literature. Finally, \cref{sec:Conclusion} concludes the work and gives a short outlook on future research directions.

        \subsection{Notation} \label{subsec:Notation}
            Unless otherwise stated, matrices are denoted by capital bold Latin or Greek letters (e.g.\ $\myvec{X}, \myvec{U}, \myvec{V}, \myvec{\Psi}, \dots$). The $i$-th row and $j$-th column of $\myvec{U}$ are written as $\myvec{U}_{i, \bullet}$ and $\myvec{U}_{\bullet,j}$ in a bold format respectively, whereas the entry in the $i$-th row and $j$-th column is indicated as $U_{ij}$ in a standard format. Vectors in $\R^n$ are typically written as small bold Latin or Greek letters (e.g.\ $\myvec{x}, \myvec{u}, \myvec{v}, \myvec{\theta}, \dots$).

            Moreover, several norms of matrices and vectors are used throughout this work. The Frobenius norm of a matrix $\myvec{M}$ and the standard Euclidean norm of a vector $\myvec{M}_{\bullet,j}$ is given by $\Vert \myvec{M}\Vert_F$ and $\Vert \myvec{M}_{\bullet,j} \Vert_2$ respectively. Furthermore, we use the notion of the $\ell_1$ norm of a matrix $\myvec{M},$ which is defined as $\Vert \myvec{M}\Vert_1 \coloneqq \sum_{i,j} \vert M_{ij} \vert.$ Furthermore, we use the notion of nonnegative matrices and write $\myvec{U}\geq 0$ or $\myvec{U}\in \R_{\geq 0}^{m\times n} $ by defining $\R_{\geq 0} \coloneqq \{ x\in \R \ \vert \ x \geq 0\}$ to indicate that $\myvec{U}$ only contains nonnegative entries. Finally, we write $\langle \myvec{x}, \myvec{y} \rangle_n$ for the standard inner product in $\R^n.$

    \section{Background and Related Work}\label{sec:Background and Related Work}
        The main task of a clustering method is to partition a given set of objects into groups, such that objects within a group are more similar to each other than objects from different groups. In a mathematical framework, this can be formulated as follows: Partition a given index set $\{1,\dots, M\}$ of a corresponding dataset $\{x_m\in \tilde{\mathcal{X}} \ \vert \ m=1,\dots,M \}$ into $K$ disjoint sets $\mathcal{I}_k \subset \{1,\dots, M\},$ such that $\cup_{k=1,\dots,K} \mathcal{I}_k = \{1,\dots, M\},$ where $\tilde{\mathcal{X}}$ is a non-empty set. Note that in this formulation, the property of similarity is still needed to be defined properly.

        $K$-means clustering ranks upon the best known partitional clustering algorithms and is also the focus of this work.
        Following \cite{Chen:2009:K-MeansCoresets}, the formal definition of $K$-means clustering can be formulated as follows: Let $(\tilde{\mathcal{X}}, \dist)$ be a metric space with $\dist(\cdot,\cdot)$ being a distance function over $\tilde{\mathcal{X}}.$ Furthermore, we define $\dist(\mathcal{X},y) \coloneqq \min_{x\in \mathcal{X}} \dist(x,y)$ for $\mathcal{X} \subset \tilde{\mathcal{X}}$ and $y\in \tilde{\mathcal{X}}.$ 
        \begin{definition}[$K$-median and $K$-means clustering]\label{def:K-median and K-means Clustering}
            A \textit{clustering} of a given set $\mathcal{X}\subset \tilde{\mathcal{X}}$ is a partition based on the set of \textit{centroids} $\mathcal{C}\coloneqq (c_1,\dots,c_K)\subset \tilde{\mathcal{X}},$ so that each point in $\mathcal{X}$ is associated to its nearest centroid $c_k.$ We assume that each point $x\in \mathcal{X}$ is associated to a weight $\weight(x)\in \R.$ The cost of $K$-median clustering and $K$-means clustering of $\mathcal{X}$ by $\mathcal{C}$ is defined as $\psi(\mathcal{C},\mathcal{X}) \coloneqq \sum_{x\in \mathcal{X}} \weight(x) \dist(\mathcal{C}, x)$ and $\omega(\mathcal{C},\mathcal{X}) \coloneqq \sum_{x\in \mathcal{X}} \weight(x) \dist(\mathcal{C}, x)^2$ respectively. The corresponding metric $K$-median (resp.\ $K$-means) problem is to find a set of centroids $\mathcal{C}\subset \tilde{\mathcal{X}},$ such that the cost $\psi(\mathcal{C},\mathcal{X})$ and $\omega(\mathcal{C},\mathcal{X})$ is minimized.
        \end{definition}
        In this setting, the so-called \textit{clusters} of the clustering are the disjoint sets $\mathcal{I}_k$ consisting of the points $x\in \mathcal{X}$ which are associated to the corresponding centroid $c_k.$ Different from \cite{Chen:2009:K-MeansCoresets}, we do not constrain in \cref{def:K-median and K-means Clustering} the set of centroids to be in $\mathcal{X}.$ In the remainder of this work, we will assume that the set $\tilde{\mathcal{X}}$ is defined by $\tilde{\mathcal{X}} \coloneqq \R^N$ along with a given dataset $\mathcal{X} \coloneqq \{\myvec{x}_m\in \R^N \ \vert \ m=1,\dots,M \}.$
        
        Furthermore, we note that according to \cref{def:K-median and K-means Clustering}, the distance function $\dist$ is a metric of the metric space $\tilde{\mathcal{X}}.$ However, $K$-means clustering algorithms oftentimes are also used with similarity measures such as the cosine similarity, which do not satisfy the properties of a metric. Moreover, it will turn out that the obtained generalized distance measures in this work will also not satisfy the usual properties of a metric (see \cref{sec:Regularized ONMF Model and Generalized K-means} for more details). However, we will refer to it as a distance measure $\dist(\cdot, \cdot)$ in the remainder of this work. 

        Thus, $K$-means clustering consists of two basic ingredients: the distance function $\dist(\cdot,\cdot),$ which gives the needed similarity property of the clustering method, as well as the centroids in $\mathcal{C}.$ The standard choice of the distance function for the classical $K$-means problem is the Euclidean distance $\dist(x,y) \coloneqq \Vert x-y\Vert_2$ leading to the typical minimization of the within-cluster variances given by
        \begin{equation}\label{eq:Background:Classical K-means Objective Function}
            \min_{\mathcal{I}_1,\dots,\mathcal{I}_K} \sum_{k = 1}^K \sum_{m\in \mathcal{I}_k} \Vert \myvec{x}_m - \myvec{c}_k \Vert_2^2.
        \end{equation}
        It can be shown, that this problem is NP-hard \cite{MahajanEtal:2012:K-MeansNPHard}. Hence, heuristic approaches are commonly used to find approximate solutions. The most widely known method is the $K$-means algorithm and is based on an alternating minimization. In a first step, after a suitable initialization of the centroids $\mathcal{C},$ the data points $\myvec{x}_m$ are associated to the nearest centroid $c_k.$ Afterwards, in the case of the Euclidean distance, the centroids $\myvec{c}_k$ are recomputed based on the mean of the associated points $\myvec{x}_m$ for $m\in \mathcal{I}_k.$ This process is repeated until the cluster assignments do not change anymore.

        In order to see the connection to matrix factorization problems, it is needed to formulate the objective function in \cref{eq:Background:Classical K-means Objective Function} in a vectorized form. To do so, we write the data points $\myvec{x}_m$ row-wise into a data matrix $\myvec{X}\in \R^{M\times N}$ so that $\myvec{X}\coloneqq [\myvec{x}_1,\dots, \myvec{x}_M]^\intercal.$ Moreover, we define the so-called \textit{cluster membership matrix} $\myvec{B}\in \{0,1\}^{M\times K},$ which is given by
        \begin{equation*}
            B_{mk} \coloneqq  \begin{cases}
                                  0 & \text{if} \ \ m\not\in \mathcal{I}_k, \\
                                  1 & \text{if} \ \ m\in \mathcal{I}_k.
                              \end{cases}
        \end{equation*}
        Note that since the sets $\mathcal{I}_k$ are disjoint to each other, each row of $\myvec{B}$ has exactly one non-zero element and gives the needed interpretation of a hard clustering. Hence, the columns of $\myvec{B}$ are orthogonal to each other, i.e.\ it holds that $\langle \myvec{B}_{\bullet,k}, \myvec{B}_{\bullet,\ell} \rangle_2 = 0$ for $k\neq \ell.$ Furthermore, note that $\vert \mathcal{I}_k \vert = \Vert \myvec{B}_{\bullet,k} \Vert_1 = \Vert \myvec{B}_{\bullet,k} \Vert_2^2.$ In addition, we introduce the diagonal matrix $\myvec{D} \coloneqq \diag(1/\vert \mathcal{I}_1 \vert, \dots, 1/\vert \mathcal{I}_K \vert)\in \R^{K\times K}.$ By using these matrices, we see that $\myvec{DB}^\intercal\myvec{X}\in \R^{K\times N}$ yields a matrix which has the centroids $c_k$ arranged in its rows. In this way, it is possible to rewrite the objective function in \cref{eq:Background:Classical K-means Objective Function} as $\Vert \myvec{X} - \myvec{BDB}^\intercal\myvec{X} \Vert_F^2$ leading to the minimization problem
        \begin{equation}\label{eq:Background:Classical K-means Objective Function as Factorization}
            \min_{\substack{\myvec{B}\in \{0,1\}^{M\times K},\\ (\myvec{BD}^{\nicefrac{1}{2}})^\intercal \myvec{BD}^{\nicefrac{1}{2}} = \myvec{I}_{K\times K}}} \Vert \myvec{X} - \myvec{BDB}^\intercal\myvec{X} \Vert_F^2,
        \end{equation}
        where $D^{\nicefrac{1}{2}}$ is defined component-wise by $D_{k\ell}^{\nicefrac{1}{2}} \coloneqq \sqrt{D_{k\ell}}$ and with $\myvec{I}_{K\times K}$ being the identity matrix of size $K\times K.$
        
        A trivial solution to this problem could be to choose $K=M$ and $\myvec{B}=\myvec{I}_{M\times M},$ which would correspond to assign each observation to its own cluster. However, this is obviously not the aim of a clustering method so that usually $K\ll \min\{M,N\}$ is chosen.

        The formulation in \cref{eq:Background:Classical K-means Objective Function as Factorization} makes the relationship between the clustering problem and a matrix factorization problem more clear. By omitting the constraints in the minimization problem \cref{eq:Background:Classical K-means Objective Function as Factorization} and additionally assuming that a nonnegative data matrix $\myvec{X}\geq 0$ is given, we automatically have that $\myvec{DB}^\intercal\myvec{X}\geq 0.$ This gives rise to the so-called Nonnegative Matrix Factorization (NMF) problem $\myvec{X} \approx \myvec{UV}$ of the data matrix, so that the factorization matrices $\myvec{U}\in \R_{\geq 0}^{M\times K}$ and $\myvec{V} \in \R_{\geq 0}^{K\times N}$ of the NMF can be compared to the cluster membership matrix $\myvec{B}$ and the so-called \textit{centroid} matrix $\myvec{DB}^\intercal\myvec{X}$ respectively. In the following, we give a definition of the general NMF problem.
        \begin{definition}[NMF]\label{def:Background:NMF}
            For a given data matrix $\myvec{X}\in \R_{\geq 0}^{M\times N}$ and a factorization rank $K\ll \min\{M,N\},$ the aim is to find two matrices $\myvec{U}\in \R_{\geq 0}^{M\times K}$ and $\myvec{V} \in \R_{\geq 0}^{K\times N},$ such that
            \begin{equation}\label{eq:def:Background:NMF}
                \myvec{X} \approx \myvec{UV} = \sum_{k=1}^K \myvec{U}_{\bullet, k}V_{k,\bullet}.
            \end{equation}
        \end{definition}
        NMF was originally introduced by Paatero and Tapper \cite{PaateroTapper:1994:PositiveMF} in 1994 as positive matrix factorization. Different from the widely-known Principal Component Analysis (PCA), NMF constraints the factorization matrices to be nonnegative. This allows a parts-based representation of the whole dataset, since each row $\myvec{X}_{m,\bullet}$ and column $\myvec{X}_{\bullet, n}$ can be represented as a superposition of the few basis vectors $\myvec{V}_{k,\bullet}$ and $\myvec{U}_{\bullet, k}$ so that $X_{\bullet,n } \approx \sum_k V_{kn} U_{\bullet, k} $ and  $X_{m,\bullet } \approx \sum_k U_{mk} V_{k, \bullet}.$ This property makes the NMF the ideal tool for nonnegative data, since the interpretability of the factorization matrices is ensured due to the additional nonnegativity constraint. NMF has been extensively used as a feature extraction and data representation tool \cite{Fernsel:2018:Survey,Leuschner:2018:NMFSupervised} as well as for clustering \cite{Fernsel:2021:ONMFTV}, compression \cite{Zhijian:2005:NMFForCompression} or even for solving inverse problems, where the NMF can be used as a joint reconstruction and feature extraction method \cite{Arridge:2020:NMFForIP}. Possible application fields include document clustering \cite{KimPark:2008:SparseNMF,PanNg:2018:CompMethod}, medical imaging \cite{Fernsel:2018:Survey,Leuschner:2018:NMFSupervised,Fernsel:2021:ONMFTV,Arridge:2020:NMFForIP}, hyperspectral unmixing \cite{HeEtal:2017:HypUnmix,FengEtal:2018:HypUnmix,FengEtal:2019:HypUnmix} and music analysis \cite{Fevotte:2009:NMFforMusic} to name just a few.

        The typical approach to find an approximate solution for the NMF is based on a variational formulation of the problem. Thus, the NMF is reformulated as a minimization problem with a suitable discrepancy term $D(\cdot, \cdot),$ which is typically chosen according to the noise distribution of the data. Furthermore, NMF problems are usually ill-posed due to the non-uniqueness of the solution \cite{Klingenberg:2009:NMFIllposedness,PhamLachmund:2020:NMFIllPosedness}. Hence, suitable regularization terms $R_j(\cdot)$ are typically added to the NMF cost function to tackle the ill-posedness of the problem and to enforce additional properties of the factorization matrices. Hence, the general minimization problem of the NMF can be written as
        \begin{equation} \label{eq:Background:NMFProblemGeneral}
            \min_{\myvec{U}, \myvec{V}\geq 0} D(\myvec{X}, \myvec{UV}) + \sum_{j=1}^J \alpha_j R_j(\myvec{U}, \myvec{V}) \eqqcolon \min_{\myvec{U}, \myvec{V}\geq 0} F(\myvec{U}, \myvec{V}),
        \end{equation}
        where $\alpha_j\geq 0$ are the regularization parameters, which control the influence of the penalty terms $R_j(\cdot).$ Typical choices for discrepancy terms are the Frobenius norm in the case of Gaussian noise, the Kullback-Leibler divergence for Poisson noise, the $\ell_1$ norm or other divergences. Regarding the penalty terms, common choices are the $\ell_2$ and $\ell_1$ regularization, which is also used in this work (see \cref{sec:Regularized ONMF Model and Generalized K-means}). Further possibilities are more problem specific and include total variation regularization and terms which enforce orthogonality of the factorization matrices or even allow a supervised classification framework \cite{Fernsel:2018:Survey,Leuschner:2018:NMFSupervised,Arridge:2020:NMFForIP,Fernsel:2021:ONMFTV}.

        For usual choices of $D$ and $R_j,$ the corresponding objective function $F$ is convex in each of the variables $\myvec{U}$ and $\myvec{V}$ but non-convex in $(\myvec{U},\myvec{V}).$ This motivates to consider alternating minimization schemes similar to the $K$-means algorithm discussed above leading to the update rules
        \begin{align}
            \myvec{U}^{[i+1]} &= \arg\min_{\myvec{U}\geq 0} F(\myvec{U},\myvec{V}^{[i]}), \label{eq:Background:AlgorithmsForNMF:AlternatingMin:U} \\
            \myvec{V}^{[i+1]} &= \arg\min_{\myvec{V}\geq 0} F(\myvec{U}^{[i+1]},\myvec{V}). \label{eq:Background:AlgorithmsForNMF:AlternatingMin:V}
        \end{align}
        For a review on the algorithm development of multiplicative updates for a variety of discrepancy and regularization terms, we refer the reader to \cite{Fernsel:2018:Survey}.

        Regarding the clustering capability of NMF and comparing $\myvec{U}$ to the cluster membership matrix $\myvec{B}$ of the $K$-means approach, the classical NMF problem with only the nonnegativity constraint on the matrices does not provide the needed hard clustering interpretability on $\myvec{U},$ since the matrix can contain multiple non-zero entries each of its rows. One typical approach to ensure this property on $\myvec{U}$ is to additionally require the matrix to be column-wise orthogonal by adding the hard constraint
        \begin{equation}\label{eq:Background:ONMF:Constraint Orthogonality:Unnormalized}
            \langle \myvec{U}_{\bullet,k}, \myvec{U}_{\bullet,\ell} \rangle_2 = 0 \quad \text{for}\ k\neq\ell.
        \end{equation}
        This leads to the problem of \textit{Orthogonal NMF} (ONMF). Occasionally, further constraints like the normalization of the columns $\myvec{U}_{\bullet,k}$ leading to
        \begin{equation}\label{eq:Background:ONMF:Constraint Orthogonality:Normalized}
            \myvec{U}^\intercal\myvec{U}=\myvec{I}_{K\times K}
        \end{equation}
        can be enforced, which will be also discussed in this work. These constraints indeed yields the desired interpretability of $\myvec{U}$ as a cluster membership matrix, since the nonnegativity constraint together with \cref{eq:Background:ONMF:Constraint Orthogonality:Unnormalized} ensure that every row of $\myvec{U}_{m,\bullet}$ only contains at most one non-zero entry $U_{m,k^*}>0$ indicating the association of the data point $\myvec{X}_{m,\bullet}$ to the cluster $\mathcal{I}_{k^*}.$ Thus, the clusters can also be written as
        \begin{equation}\label{eq:Background:Ik via U}
            \mathcal{I}_k = \{m\in \{1,\dots,M\} \ \vert \ U_{mk}>0 \}.
        \end{equation}
        In this setting, the matrix $\myvec{V}$ can then be interpreted as the centroid matrix, which contains the centroids in its rows.

        Throughout the literature, many relationships between different kinds of $K$-means and ONMF models could be shown. One of the pioneering works is the one by Ding et al.\ in \cite{Ding:2005:NMFKMeans}, which describes equivalences between Kernel $K$-means and symmetric NMF $\myvec{X}\approx \myvec{UU}^\intercal,$ bipartite graph $K$-means clustering and the bi-orthogonal NMF problem given by the minimization problem
        \begin{equation*}
            \min_{\substack{U,V\geq 0}} \Vert \myvec{X} - \myvec{UV} \Vert_F^2, \quad \text{s.t.} \ \myvec{U}^\intercal\myvec{U}=\myvec{I}_{K\times K},\ \ \myvec{V}\myvec{V}^\intercal=\myvec{I}_{K\times K},
        \end{equation*}
        as well as the classical $K$-means clustering and ONMF with the stronger constraint on $\myvec{U}$ by considering $U_{mk}=1/\sqrt{\vert \mathcal{I}_k\vert}$ for $m\in \mathcal{I}_k$ and $U_{mk}=0$ otherwise. Numerous other works followed with relationships between Nonnegative Matrix Tri-Factorizations and simultaneous row and column clustering approaches with applications to document clustering problems as well as connections between relaxed $K$-means clustering models and semi-NMF, convex NMF and Kernel NMF \cite{LiDing:2006:Relationships,Ding:2006:triFac}. Furthermore, \cite{Pompili:2014:CompMethod} shows the equivalence between a weighted variant of spherical $K$-means and the ONMF model with \cref{eq:Background:ONMF:Constraint Orthogonality:Normalized} as a hard constraint on $\myvec{U}.$ Further results include the clustering interpretability of sparse NMF \cite{KimPark:2008:SparseNMF} and the relation of projective NMF to $K$-means clustering \cite{Yuan:2009:ProjectiveNMF}.

        Besides of the theoretical interest to study these relationships, they also have some practical relevance since NMF models can have several advantages over classical $K$-means clustering methods. For instance, NMF models can do both hard as well as soft clustering and are able to perform a clustering of the rows and columns simultaneously. For more information on the clustering capabilities of NMF approaches, the relationships to $K$-means clustering and the development of algorithms, we refer the interested reader to both survey articles \cite{LiDing:2014:NMFClusteringSurvey,Turkmen:2015:NMFClusteringSurvey}.

        However, these works do not consider any regularization terms in their NMF models, which is the focus of this work. Furthermore, the derivation of these relationships in the works throughout the literature are mostly based on the comparison of the objective functions and the constraints of the considered $K$-means and NMF model. This is in contrast to the approach used in this work, where the focus lies on obtaining connections between regularized ONMF and generalized $K$-means models by directly deriving the distance measures and centroids of the $K$-means approach based on the considered ONMF model. The used derivation framework exploits the typical alternating minimization scheme of NMF algorithms in \cref{eq:Background:AlgorithmsForNMF:AlternatingMin:U} and \cref{eq:Background:AlgorithmsForNMF:AlternatingMin:V} and uses first-order conditions of the ONMF objective function.
        This technique offers a significantly simpler method to derive connections between $K$-means and NMF compared to the ones used throughout the literature. Furthermore, it can be generalized to regularized ONMF models and is able to directly extract the distance measure and centroids of the $K$-means approach. In \cref{sec:Regularized ONMF Model and Generalized K-means}, we consider a regularized ONMF model with an elastic net regularization on both matrices $\myvec{U}$ and $\myvec{V}.$ Hence, from a viewpoint of the regularization theory in inverse problems, the used framework also allows to see the effect of the $\ell_1$ and $\ell_2$ penalty terms on the obtained distance measures and the centroids of the $K$-means method.   

        Another motivational aspect concerns the distance measure for $K$-means clustering. It is well-known that an appropriate choice of the distance measure depending on the considered application is vital for the performance of the clustering algorithm. However, it is also known that $K$-means clustering algorithms typically suffer from the so-called uniformization effect, i.e.\ the algorithm tends to produce clusters with relatively balanced sizes. The work \cite{WuEtal:2007:K-meansDistanceUniformization} could show that for a so-called $K$-means distance, which is a generalization of the Bregman divergence, the usual $K$-means algorithm suffer from the uniformization effect. However, the presented framework in this paper also allows to derive generalized $K$-means models with distance measures, which do not belong to the family of Bregman divergences and could lead to $K$-means approaches, which do not suffer from this negative effect.

    \section{Regularized ONMF and \texorpdfstring{$K$}{K}-means Clustering} \label{sec:Regularized ONMF Model and Generalized K-means}
        In this section, we introduce the considered ONMF problem and study its relation to generalized $K$-means models by deriving the corresponding distance measures and centroids. Regarding the ONMF model, we consider the $\ell_1$ norm and the $\ell_2$ norm for the discrepancy term and an elastic net regularization on $\myvec{U}$ and $\myvec{V}$ in each case. Hence, the objective functions of the considered ONMF models are
        \begin{align}
            F_1(\myvec{U}, \myvec{V}) &\coloneqq \underbrace{\Vert \myvec{X} - \myvec{UV} \Vert_1}_{\eqqcolon D_1(\myvec{X},\myvec{UV})} + \underbrace{\lambda_{\myvec{U}} \Vert \myvec{U} \Vert_1 + \mu_{\myvec{U}} \Vert \myvec{U}\Vert_F^2 + \lambda_{\myvec{V}} \Vert \myvec{V} \Vert_1 + \mu_{\myvec{V}} \Vert \myvec{V}\Vert_F^2}_{\eqqcolon R(\myvec{U},\myvec{V})},\label{eq:ONMF Model:Cost Function F1 for l1}\\
            F_2(\myvec{U}, \myvec{V}) &\coloneqq \underbrace{\Vert \myvec{X} - \myvec{UV} \Vert_F^2}_{\eqqcolon D_2(\myvec{X},\myvec{UV})} + \lambda_{\myvec{U}} \Vert \myvec{U} \Vert_1 + \mu_{\myvec{U}} \Vert \myvec{U}\Vert_F^2 + \lambda_{\myvec{V}} \Vert \myvec{V} \Vert_1 + \mu_{\myvec{V}} \Vert \myvec{V}\Vert_F^2,\label{eq:ONMF Model:Cost Function F2 for l2}
        \end{align}
        where $\myvec{X}\coloneqq [\myvec{x}_1,\dots, \myvec{x}_M]^\intercal \in \R_{\geq 0}^{M\times N}$ is a given data matrix with data points $\myvec{x}_m\in \R_{\geq 0}^{N},$ so that $\myvec{X}_{m,\bullet} = (\myvec{x}_m)^\intercal.$ Furthermore, $\myvec{U} \coloneqq [\myvec{u}_1,\dots, \myvec{u}_K] \in \R_{\geq 0}^{M\times K}$ and $\myvec{V} \coloneqq [\myvec{v}_1,\dots, \myvec{v}_K]^\intercal \in \R_{\geq 0}^{K\times N}$ are the factorization matrices of the ONMF problem with the short notations $\myvec{u}_k = \myvec{U}_{\bullet,k}$ and $\myvec{v}_k = (\myvec{V}_{k,\bullet})^\intercal.$ Moreover, $\lambda_{\myvec{U}}, \lambda_{\myvec{V}}, \mu_{\myvec{U}}, \mu_{\myvec{V}}\geq 0$ are the regularization parameters and control the influence of the corresponding penalty terms in the objective function. In the most general case, we consider the ONMF minimization problem
        \begin{equation}\label{eq:ONMF Problem}
            \begin{aligned}
                \min_{\myvec{U},\myvec{V}\geq 0}\quad &F_i(\myvec{U},\myvec{V})\\
                \text{s.t.}\quad &\cref{eq:ONMF Model:Constraint U}\ \text{holds}
            \end{aligned}
        \end{equation}
        for $i\in \{1,2\}$ with the additional constraint
        \begin{equation}\label{eq:ONMF Model:Constraint U}
            \langle \myvec{u}_k,\myvec{u}_\ell \rangle_2 = 0 \quad \text{for} \quad k\neq l \tag{C1}
        \end{equation}
        on $\myvec{U}$ to ensure the needed clustering interpretability as a cluster membership matrix as described in \cref{sec:Background and Related Work}. Occasionally, we will constrain the matrix $\myvec{U}$ further than in \cref{eq:ONMF Model:Constraint U} to discuss special cases (see \cref{subsec:Alternating Minimization and Separability} and \cref{subsec:Evaluation and Special Cases}). 
        \subsection{Alternating Minimization and Separability} \label{subsec:Alternating Minimization and Separability}
            In this section, we describe the basic strategy to derive the connections between the ONMF and $K$-means models and introduce some basic tools used in the following Sections.

            As described in \cref{sec:Background and Related Work}, the general framework used in this work is based on the alternating minimization scheme showed in \cref{eq:Background:AlgorithmsForNMF:AlternatingMin:U} and \cref{eq:Background:AlgorithmsForNMF:AlternatingMin:V}. As the matrix $\myvec{V}$ can be compared to a centroid matrix in the setting of the considered ONMF model (see \cref{sec:Background and Related Work}), the minimization of the ONMF problem in \cref{eq:ONMF Problem} with respect to $\myvec{V}$ for fixed $\myvec{U}$ leads to the characterization of the centroid of the $K$-means approach. Furthermore, it will turn out that the distance measure of the $K$-means method can be identified via the minimization of \cref{eq:ONMF Problem} with respect to $\myvec{U}$ for fixed $\myvec{V}.$ However, in order be able to compute the solutions of both minimization problems, the objective functions $F_i$ need to satisfy specific separability properties.

            In the case of the minimization with respect to $\myvec{V},$ it is easy to see that both $F_i$ can be written as
            \begin{equation*}
                F_i(\myvec{U}, \myvec{V}) = \sum_{k=1}^K D_i(\myvec{X}_k, \tilde{\myvec{u}}_k \myvec{v}_k^\intercal) + R(\myvec{u}_k, \myvec{v}_k),
            \end{equation*}
            where $ \myvec{X}_k\in \R_{\geq 0}^{\vert \mathcal{I}_k \vert \times N}$ denotes the submatrix containing the rows of $ \myvec{X} $ with the indices in $\mathcal{I}_k$ given by the relationship in \cref{eq:Background:Ik via U}. Furthermore, $\tilde{\myvec{u}}_k \in \R_{\geq 0}^{\vert \mathcal{I}_k \vert}$ denotes the corresponding reduced column vector of $ \myvec{u}_k $ of non-zero entries. Hence, the minimization of the $F_i$ with respect to $ \myvec{V} $ separates with the rows $ \myvec{v}_k. $ Therefore, the centroids $\myvec{v}_k$ of the $K$-means approach can be computed by solving the minimization problems
            \begin{equation}\label{eq:AlternatingMin:V:Vectorized}
                \min_{\myvec{v}_k \geq 0} D_i(\myvec{X}_k, \tilde{\myvec{u}}_k \myvec{v}_k^\intercal) + \mu_{\myvec{V}} \Vert \myvec{v}_k\Vert_2^2 + \lambda_{\myvec{V}} \Vert \myvec{v}_k \Vert_1,
            \end{equation}
            since the terms independent from $\myvec{V}$ can be omitted. We can further reduce the minimization problem to a set of scalar minimization problems, since
            \begin{equation*}
                D_i(\myvec{X}_k, \tilde{\myvec{u}}_k \myvec{v}_k^\intercal) + R(\myvec{u}_k, \myvec{v}_k) = \sum_{n=1}^N \sum_{m\in \mathcal{I}_k} D_i(X_{mn}, U_{mk}V_{kn}) + \mu_{\myvec{V}} V_{kn}^2 + \lambda_{\myvec{V}} \vert V_{kn} \vert.
            \end{equation*}
            This leads to the scalar minimization problems
            \begin{equation}\label{eq:AlternatingMin:V:Scalar}
                \min_{V_{kn} \geq 0} \sum_{m\in \mathcal{I}_k} D_i(X_{mn}, U_{mk}V_{kn}) + \mu_{\myvec{V}} V_{kn}^2 + \lambda_{\myvec{V}} \vert V_{kn} \vert \eqqcolon \min_{V_{kn} \geq 0} \varphi_{kn}(V_{kn})
            \end{equation}
            with the objective functions
            \begin{equation}\label{eq:AlternatingMin:V:Scalar:ObjectiveFunction}
                \varphi_{kn}(t) = \sum_{m\in \mathcal{I}_k} D_i(X_{mn}, U_{mk}t) + \lambda_{\myvec{V}} \vert t \vert + \mu_{\myvec{V}} t^2,
            \end{equation}
            which are considered in the following Sections \ref{subsec:Classical Discrepancy} and \ref{subsec:l1 Discrepancy}.

            Regarding the minimization task with respect to $\myvec{U}$ for fixed $\myvec{V},$ the problem also separates for the individual entries of $\myvec{U}.$ Due to the constraint in \cref{eq:ONMF Model:Constraint U} in combination with the nonnegativity of $\myvec{U},$ each row of $\myvec{U}$ has at most one positive entry. Thus, the $m$-th row of $\myvec{UV}$ is a multiple of one of the rows of $\myvec{V}$ and is compared in the discrepancy term with $\myvec{X}_{m,\bullet}.$ Hence, the objective functions $F_i$ can be written as
            \begin{equation*}
                F_i(\myvec{U}, \myvec{V}) = \sum_{m=1}^M \left[ D_i(\myvec{X}_{m,\bullet}, U_{m,\pi(m)}\myvec{V}_{\pi(m),\bullet}) + \lambda_{\myvec{U}} \vert U_{m,\pi(m)} \vert + \mu_{\myvec{U}} U_{m,\pi(m)}^2 \right] + \lambda_{\myvec{V}} \Vert \myvec{V} \Vert_1 + \mu_{\myvec{V}} \Vert \myvec{V}\Vert_F^2,
            \end{equation*}
            where $\pi(m)\in \{1,\dots, K\}$ denotes the index of the non-zero entry in the $m$-th row of $\myvec{U}$ assuming that in each row of $\myvec{U},$ there exists exactly one non-zero entry $U_{m,\pi(m)}>0.$ This motivates to consider the minimization problems
            \begin{equation}\label{eq:AlternatingMin:U:Scalar}
                \min_{U_{m,\pi(m)}>0} D_i(\myvec{X}_{m,\bullet}, U_{m,\pi(m)}\myvec{V}_{\pi(m),\bullet}) + \mu_{\myvec{U}} U_{m,\pi(m)}^2 + \lambda_{\myvec{U}} \vert U_{m,\pi(m)} \vert \eqqcolon \min_{U_{m,\pi(m)}>0} \psi_{m,\pi(m)}(U_{m,\pi(m)})
            \end{equation}
            with the objective functions
            \begin{equation}\label{eq:AlternatingMin:U:Scalar:ObjectiveFunction}
                \psi_{m,\pi(m)}(t) = D_i(\myvec{X}_{m,\bullet}, t\myvec{V}_{\pi(m),\bullet}) + \mu_{\myvec{U}} t^2 + \lambda_{\myvec{U}} \vert t \vert
            \end{equation}
            omitting the terms independent from $U_{m,\pi(m)}.$ This minimization problem characterizes the distance measure of the $K$-means approach and yields
            \begin{equation}\label{eq:AlternatingMin:U:DistanceFunction}
                \dist(\myvec{x}_m, \myvec{v}_{\pi(m)}) = \left( \min_{t> 0} \psi_{m,\pi(m)}(t) \right)^{\nicefrac{1}{i}},
            \end{equation}
            where the index $i$ specifies the considered objective function $F_i.$ The index $\pi(m),$ which indicates the non-zero entry $U_{m,\pi(m)},$ can be computed by comparing the $m$-th data vector $\myvec{x}_m$ with all centroids of $\myvec{v}_k$ leading to
            \begin{equation*}
                \pi(m) = \arg\min_{k\in \{1,\dots,K\}} \dist(\myvec{x}_m, \myvec{v}_k).
            \end{equation*}
            The value of the non-zero entry in the $m$-th row of the matrix $ \myvec{U} $ is finally given by
            \begin{equation*}
                U_{m,\pi(m)} = \arg \min_{t>0} \psi_{m,\pi(m)}(t).
            \end{equation*}

            We note that both scalar minimization problems in \cref{eq:AlternatingMin:V:Scalar} and \cref{eq:AlternatingMin:U:Scalar} are strictly convex for $\mu_{\myvec{V}}>0$ and $\mu_{\myvec{U}}>0$ respectively due to the quadratic penalty term of the elastic net regularization.

            The above described relationships between the ONMF model and the distance function as well as the centroid of the $K$-means approach can be seen more intuitively by considering the additional constraint $i=2,\ U_{mk}\in \{0,1\}$ and $\lambda_{\myvec{U}} = \lambda_{\myvec{V}} = \mu_{\myvec{U}} = \mu_{\myvec{V}} = 0$ as a special case of the ONMF model. For the minimization with respect to $\myvec{V},$ this yields the minimization problem $\min_{\myvec{v}_k\geq 0} \Vert \myvec{X}_k - \tilde{\myvec{u}}_k \myvec{v}_k^\intercal\Vert_F^2.$ Using the zero gradient condition, this leads to $\myvec{v}_k = \sum_{m\in \mathcal{I}_k}\myvec{x}_m / \vert \mathcal{I}_k \vert,$ i.e.\ the mean of the data points of the index set given in $\mathcal{I}_k.$ For the minimization problem with respect to $\myvec{U},$ we consider \cref{eq:AlternatingMin:U:DistanceFunction} with the additional constraint $t\in \{0,1\}$ leading to the usual Euclidean distance, i.e.\ $\dist(\myvec{x}_m, \myvec{v}_k) = \Vert\myvec{x}_m - \myvec{v}_k\Vert_2.$ Hence, this ONMF model leads to the classical $K$-means approach with the Euclidean distance function and mean centroids.

            For the subsequent analysis, we introduce the notion of a generalized median function. To do so, we consider the vectors $\myvec{v,w}\in \R^{N}$ and $\lambda, \mu\geq 0$ along with the continuous and convex function
            \begin{equation}\label{eq:GeneralizedMedian:FunctionF}
                f(t) = \sum_{n=1}^N \vert v_n - w_n t \vert + \mu t^2 + \lambda \vert t\vert.
            \end{equation}
            It is easy to see that minimizers of $f$ exist and that the set of minimizers $\arg\min_t f(t)$ either contains a unique element or is a closed interval. For both cases, we define the weighted, regularized median
            \begin{equation}\label{eq:GeneralizedMedian}
                \med_{\myvec{w}}^{\ell_1(\lambda), \ell_2(\mu)} (\myvec{v}) \in \arg\min_t f(t) \subset \R
            \end{equation}
            as the midpoint of the set of minimizers given by $\arg\min_t f(t).$ Note that for $\lambda = \mu = 0$ and $w_n = 1$ for all $n\in \{1,\dots,N\},$ the weighted regularized median is the classical median, so that $\med_{\myvec{w}}^{\ell_1(\lambda), \ell_2(\mu)} (\myvec{v}) = \med(\myvec{v}).$
        \subsection{Classical Discrepancy Term} \label{subsec:Classical Discrepancy}
            In this section, we derive the distance measures and the centroids based on the ONMF problem in \cref{eq:ONMF Problem} for $i=2.$ For the minimization with respect to $\myvec{V},$ we consider the scalar minimization problem
            \begin{equation*}
                \min_{t\geq 0} \sum_{m\in \mathcal{I}_k} (X_{mn} - U_{mk}t)^2 + \mu_{\myvec{V}} t^2 + \lambda_{\myvec{V}} \vert t \vert \eqqcolon \min_{t\geq 0} \varphi_{kn}(t)
            \end{equation*}
            based on the objective function in \cref{eq:AlternatingMin:V:Scalar:ObjectiveFunction}. To find a minimizer of the above problem, we examine the zero gradient condition
            \begin{equation*}
                0\in \partial \varphi_{kn}(t) = \sum_{m\in \mathcal{I}_k} 2 (U_{mk}t - X_{mn})U_{mk} + 2 \mu_{\myvec{V}} t + \lambda_{\myvec{V}} \partial \vert \cdot \vert (t),
            \end{equation*}
            where $\partial \vert \cdot \vert$ is the subdifferential of the absolute value function. This can be equivalently written as
            \begin{equation*}
                \frac{\sum_{m\in \mathcal{I}_k} X_{mn} U_{mk}}{\Vert \myvec{u}_k \Vert_2^2 + \mu_{\myvec{V}}} \in t + \frac{\lambda_{\myvec{V}}}{2(\Vert \myvec{u}_k \Vert_2^2 + \mu_{\myvec{V}})} \partial \vert \cdot \vert (t) \eqqcolon (\id + \gamma \partial \vert \cdot \vert)(t)
            \end{equation*}
            with $\gamma \coloneqq \lambda_{\myvec{V}} / (2 \Vert \myvec{u}_k \Vert_2^2 + 2 \mu_{\myvec{V}}).$ By using the soft thresholding function $\tau_\gamma \coloneqq (\id + \gamma \partial \vert \cdot \vert)^{-1},$ we finally obtain
            \begin{equation}\label{eq:l2:MinV:ArgMin}
                V_{kn} = \arg\min_{t\geq 0} \varphi_{kn}(t) = \tau_\gamma \left( \frac{\sum_{m\in \mathcal{I}_k} X_{mn} U_{mk}}{\Vert \myvec{u}_k \Vert_2^2 + \mu_{\myvec{V}}} \right) = \begin{cases}
                    \frac{\sum_{m\in \mathcal{I}_k} X_{mn} U_{mk}}{\Vert \myvec{u}_k \Vert_2^2 + \mu_{\myvec{V}}} - \gamma, & \gamma \leq \frac{\sum_{m\in \mathcal{I}_k} X_{mn} U_{mk}}{\Vert \myvec{u}_k \Vert_2^2 + \mu_{\myvec{V}}},\\
                    0 & \text{else},
                \end{cases}
            \end{equation}
            where the case
            \begin{equation*}
                \frac{\sum_{m\in \mathcal{I}_k} X_{mn} U_{mk}}{\Vert \myvec{u}_k \Vert_2^2 + \mu_{\myvec{V}}} \leq -\gamma
            \end{equation*}
            can be omitted due to the nonnegativity constraint.

            For the minimization with respect to $\myvec{U},$ we consider the minimization task
            \begin{equation}\label{eq:l2:U:MinTask}
                \min_{t>0} \sum_{n=1}^N (X_{m,n} - V_{\pi(m),n} t)^2 + \mu_{\myvec{U}} t^2 + \lambda_{\myvec{U}} \vert t \vert \eqqcolon \min_{t>0} \psi_{m,\pi(m)}(t)
            \end{equation}
            based on the objective function in \cref{eq:AlternatingMin:U:Scalar:ObjectiveFunction}. Similar as in the case before, the solution can be obtained via the zero gradient condition
            \begin{equation*}
                0\in \sum_{n=1}^N 2(V_{\pi(m),n} t - X_{m,n})V_{\pi(m),n} + 2\mu_{\myvec{U}} t + \lambda_{\myvec{U}} \partial \vert \cdot \vert (t),
            \end{equation*}
            which can be equivalently written as
            \begin{equation*}
                \frac{\langle \myvec{x}_m, \myvec{v}_{\pi(m)} \rangle_2}{\Vert \myvec{v}_{\pi(m)}\Vert_2^2 + \mu_{\myvec{U}}} \in t + \frac{\lambda_{\myvec{U}}}{2(\Vert \myvec{v}_{\pi(m)}\Vert_2^2 + \mu_{\myvec{U}})} \partial \vert \cdot \vert (t) \eqqcolon (\id + \gamma \partial \vert \cdot \vert)(t)
            \end{equation*}
            with $\gamma \coloneqq \lambda_{\myvec{U}} / (2\Vert \myvec{v}_{\pi(m)}\Vert_2^2 + 2\mu_{\myvec{U}}).$ As before, by using the soft thresholding function $\tau_\gamma,$ we finally obtain
            \begin{equation}\label{eq:l2:MinU:ArgMin}
                U_{m,\pi(m)} = \arg\min_{t>0} \psi_{m,\pi(m)}(t) = \tau_\gamma \left( \frac{\langle \myvec{x}_m, \myvec{v}_{\pi(m)} \rangle_2}{\Vert \myvec{v}_{\pi(m)}\Vert_2^2 + \mu_{\myvec{U}}} \right) = \begin{cases}
                    \frac{\langle \myvec{x}_m, \myvec{v}_{\pi(m)} \rangle_2}{\Vert \myvec{v}_{\pi(m)}\Vert_2^2 + \mu_{\myvec{U}}} - \gamma, & \gamma \leq \frac{\langle \myvec{x}_m, \myvec{v}_{\pi(m)} \rangle_2}{\Vert \myvec{v}_{\pi(m)}\Vert_2^2 + \mu_{\myvec{U}}},\\
                    0 & \text{else}.
                \end{cases}
            \end{equation}
            Based on \cref{eq:AlternatingMin:U:DistanceFunction} and by inserting \cref{eq:l2:MinU:ArgMin} into $\psi_{m,\pi(m)},$ we get the distance measure
            \begin{equation}\label{eq:l2:DistanceMeasure}
                \begin{aligned}
                    \dist(\myvec{x}_m, \myvec{v}_{\pi(m)})^2 = \left\Vert \myvec{x}_m - \tau_\gamma \left( \frac{\langle \myvec{x}_m, \myvec{v}_{\pi(m)} \rangle_2}{\Vert \myvec{v}_{\pi(m)}\Vert_2^2 + \mu_{\myvec{U}}} \right) \myvec{v}_{\pi(m)} \right\Vert_2^2 &+ \mu_{\myvec{U}} \left( \tau_\gamma \left( \frac{\langle \myvec{x}_m, \myvec{v}_{\pi(m)} \rangle_2}{\Vert \myvec{v}_{\pi(m)}\Vert_2^2 + \mu_{\myvec{U}}} \right)\right)^2\\
                    &+ \lambda_{\myvec{U}} \left\vert \tau_\gamma \left( \frac{\langle \myvec{x}_m, \myvec{v}_{\pi(m)} \rangle_2}{\Vert \myvec{v}_{\pi(m)}\Vert_2^2 + \mu_{\myvec{U}}} \right) \right\vert
                \end{aligned}
            \end{equation}
            between a data point $\myvec{x}_m$ and a centroid $\myvec{v}_{\pi(m)}.$
        \subsection{\texorpdfstring{$\ell_1$}{l1} Discrepancy Term} \label{subsec:l1 Discrepancy}
            Different from the classical Frobenius norm as a discrepancy term, the $\ell_1$ discrepancy term leads to a more robust NMF formulaton with respect to noise and outliers \cite{KongEtal:2011:RobustNMF}.

            Regarding the minimization with respect to $\myvec{V},$ we consider the scalar minimization problem
            \begin{equation*}
                \min_{t\geq 0} \sum_{m\in \mathcal{I}_k} \vert X_{mn} - U_{mk}t \vert + \mu_{\myvec{V}} t^2 + \lambda_{\myvec{V}} \vert t \vert \eqqcolon \min_{t\geq 0} \varphi_{kn}(t)
            \end{equation*}
            with the objective function given in \cref{eq:AlternatingMin:V:Scalar:ObjectiveFunction}. The latter has the same structure as in \cref{eq:GeneralizedMedian:FunctionF} and can be solved by taking the weighted and regularized median defined in \cref{eq:GeneralizedMedian}. Hence, we have that
            \begin{equation}\label{eq:l1:MinV:ArgMin}
                V_{kn} = \med_{\tilde{\myvec{U}}_{\bullet,k}}^{\ell_1(\lambda_{\myvec{V}}), \ell_2(\mu_{\myvec{V}})} (\tilde{\myvec{X}}_{\bullet,n}) \in \arg\min_{t\geq 0} \varphi_{kn}(t),
            \end{equation}
            where $\tilde{\myvec{U}}_{\bullet,k}, \tilde{\myvec{X}}_{\bullet,n} \in \R_{\geq 0}^{\vert \mathcal{I}_k\vert}$ are the corresponding reduced vectors with the entries $\tilde{U}_{\tilde{m},k} \in \{ U_{mk} \ \vert\ m\in \mathcal{I}_k\}$ and $\tilde{X}_{\tilde{m},n} \in \{ X_{mn} \ \vert\ m\in \mathcal{I}_k \}$ for $\tilde{m}=1,\dots,\vert \mathcal{I}_k \vert.$

            Similarly, for the minimization with respect to $\myvec{U},$ we consider the minimization task
            \begin{equation*}
                \min_{t>0} \sum_{n=1}^N \vert X_{m,n} - V_{\pi(m),n} t\vert + \mu_{\myvec{U}} t^2 + \lambda_{\myvec{U}} \vert t \vert \eqqcolon \min_{t>0} \psi_{m,\pi(m)}(t)
            \end{equation*}
            based on the objective function in \cref{eq:AlternatingMin:U:Scalar:ObjectiveFunction} to obtain the non-zero entry $U_{m,\pi(m)}$ in the $m$-th row of $\myvec{U}.$ This also has the same structure as in \cref{eq:GeneralizedMedian:FunctionF} and thus yields
            \begin{equation}\label{eq:l1:MinU:ArgMin}
                U_{m,\pi(m)} = \med_{\myvec{v}_{\pi(m)}}^{\ell_1(\lambda_{\myvec{U}}), \ell_2(\mu_{\myvec{U}})} (\myvec{x}_m) \in \arg\min_{t>0} \psi_{m,\pi(m)}(t),
            \end{equation}
            which can be seen as the weighted median of the entries in $\myvec{x}_m$ with elastic net regularization. Finally, based on \cref{eq:AlternatingMin:U:DistanceFunction} and by inserting \cref{eq:l1:MinU:ArgMin} into the function $\psi_{m,\pi(m)},$ we obtain the distance measure
            \begin{equation}\label{eq:l1:DistanceMeasure}
                \begin{aligned}
                    \dist(\myvec{x}_m, \myvec{v}_{\pi(m)}) = \left\Vert \myvec{x}_m - \med_{\myvec{v}_{\pi(m)}}^{\ell_1(\lambda_{\myvec{U}}), \ell_2(\mu_{\myvec{U}})} (\myvec{x}_m) \myvec{v}_{\pi(m)} \right\Vert_1 &+ \mu_{\myvec{U}} \left( \med_{\myvec{v}_{\pi(m)}}^{\ell_1(\lambda_{\myvec{U}}), \ell_2(\mu_{\myvec{U}})} (\myvec{x}_m)\right)^2\\
                    &+ \lambda_{\myvec{U}} \left\vert \med_{\myvec{v}_{\pi(m)}}^{\ell_1(\lambda_{\myvec{U}}), \ell_2(\mu_{\myvec{U}})} (\myvec{x}_m) \right\vert
                \end{aligned}
            \end{equation}
            between a data point $\myvec{x}_m$ and a centroid $\myvec{v}_{\pi(m)}.$ A more intuitive understanding and the consideration of special cases of the generalized distance measures in \cref{eq:l1:DistanceMeasure} and \cref{eq:l2:DistanceMeasure} can be found in the following \cref{subsec:Evaluation and Special Cases}.
        \subsection{Special Cases} \label{subsec:Evaluation and Special Cases}
            In this section, we give a more intuitive understanding of the obtained distance measures and centroids in the previous sections and discuss some special cases.

            First, we note that the obtained distance measures in \cref{eq:l1:DistanceMeasure} and \cref{eq:l2:DistanceMeasure} do not satisfy the usual properties of a metric in a metric space as in \cref{def:K-median and K-means Clustering}. However, the distance measure $\dist(\myvec{x}_m, \myvec{v}_{\pi(m)})$ for both the $\ell_1$ and $\ell_2$ discrepancy term can be interpreted as a regularized $\ell_1$ and $\ell_2$ projection of the data point $\myvec{x}_m$ onto $\myvec{v}_{\pi(m)}$ respectively.

            In the case of the $\ell_2$ discrepancy term, we obtain due to \cref{eq:l2:MinU:ArgMin} two different distance measures in the case of sparsity regularization with $\lambda_{\myvec{U}}>0.$ In the case of $\gamma > \langle \myvec{x}_m, \myvec{v}_{\pi(m)} \rangle_2 / (\Vert \myvec{v}_{\pi(m)}\Vert_2^2 + \mu_{\myvec{U}}),$ which means $\lambda_{\myvec{U}}/2 > \langle \myvec{x}_m, \myvec{v}_{\pi(m)} \rangle_2,$ this leads to approximately orthogonal $\myvec{x}_m$ and $\myvec{v}_{\pi(m)}$ if $\lambda_{\myvec{U}}$ is small and to $\tau_\gamma ( \langle \myvec{x}_m, \myvec{v}_{\pi(m)} \rangle_2 / (\Vert \myvec{v}_{\pi(m)}\Vert_2^2 + \mu_{\myvec{U}}) ) = 0,$ so that we obtain the distance measure $\dist(\myvec{x}_m, \myvec{v}_{\pi(m)})^2 = \Vert \myvec{x}_m \Vert_2^2.$ For the other case $\lambda_{\myvec{U}}/2 \leq \langle \myvec{x}_m, \myvec{v}_{\pi(m)} \rangle_2,$ we obtain the more complex distance measure
            \begin{align}
                \dist(\myvec{x}_m, \myvec{v}_{\pi(m)})^2 &= \left\Vert \myvec{x}_m - \left( \frac{\langle \myvec{x}_m, \myvec{v}_{\pi(m)} \rangle_2}{\Vert \myvec{v}_{\pi(m)}\Vert_2^2 + \mu_{\myvec{U}}} - \gamma \right) \myvec{v}_{\pi(m)} \right\Vert_2^2 &&\hspace{-6ex}+ \mu_{\myvec{U}} \left( \frac{\langle \myvec{x}_m, \myvec{v}_{\pi(m)} \rangle_2}{\Vert \myvec{v}_{\pi(m)}\Vert_2^2 + \mu_{\myvec{U}}} - \gamma \right)^2 \nonumber \\
                & &&\hspace{-6ex}+ \lambda_{\myvec{U}} \bigg\vert \underbrace{\frac{\langle \myvec{x}_m, \myvec{v}_{\pi(m)} \rangle_2}{\Vert \myvec{v}_{\pi(m)}\Vert_2^2 + \mu_{\myvec{U}}} - \gamma}_{\geq 0} \bigg\vert \nonumber\\
                &= \Vert \myvec{x}_m \Vert_2^2 - \frac{(\lambda_{\myvec{U}} - 2 \langle \myvec{x}_m, \myvec{v}_{\pi(m)} \rangle_2)^2}{4 (\Vert \myvec{v}_{\pi(m)}\Vert_2^2 + \mu_{\myvec{U}})}.\label{eq:l2:DistanceMeasure:SpecialCases} &&
            \end{align}
            Note that in both cases, the distance measure does not have the usual properties of a metric in a metric space as $\dist(\myvec{x}, \myvec{x}) \neq 0$ for some vector $\myvec{x}\in \R_{\geq 0}^N.$ Besides the classical Tikhonov regularization with the corresponding regularization parameters $\mu_{\myvec{U}}$ and $\mu_{\myvec{V}},$ it is interesting to examine the connection between the obtained distance measure and the sparsity regularization with $\lambda_{\myvec{U}}>0.$ In the case of a high sparsity regularization along with the case $\lambda_{\myvec{U}}/2 > \langle \myvec{x}_m, \myvec{v}_{\pi(m)} \rangle_2,$ we obtain a distance measure $\dist(\myvec{x}_m, \myvec{v}_{\pi(m)})$ which is independent from the centroid $\myvec{v}_{\pi(m)}$ so that in these cases, the corresponding data point can be assigned arbitrarily to any centroid. This can be also inferred based on the sparsity regularization on $\myvec{U}.$ Due to the nonnegativity of $\myvec{U}$ together with the orthogonality constraint in \cref{eq:ONMF Model:Constraint U}, the matrix is already sparse and only contains at most one non-zero entry $U_{m^*k^*}$ in each row indicating the association of the data point $\myvec{X}_{m^*,\bullet}$ to the cluster $\mathcal{I}_{k^*}.$ However, if $\lambda_{\myvec{U}}$ is set sufficiently large, the additional sparsity regularization can lead to rows $\myvec{U}_{m,\bullet}$ without any non-zero entry and yields the same interpretation, namely that the corresponding data point can be assigned to an arbitrary cluster.

            Before examining unregularized versions of the ONMF models, it is interesting to consider the distance measure in \cref{eq:l2:DistanceMeasure:SpecialCases} with $\mu_{\myvec{U}} = 0.$ It can be easily shown that this leads to
            \begin{equation*}
                \dist(\myvec{x}_m, \myvec{v}_{\pi(m)})^2 = \Vert \myvec{x}_m \Vert_2^2 \sin^2 \sphericalangle(\myvec{x}_m, \myvec{v}_{\pi(m)}) + \frac{\lambda_{\myvec{U}}}{\Vert \myvec{v}_{\pi(m)} \Vert_2^2} \left( \langle \myvec{x}_m, \myvec{v}_{\pi(m)} \rangle_2 - \frac{\lambda_{\myvec{U}}}{4} \right),
            \end{equation*}
            which explicitly shows the angle dependence of the distance measure. Here, we use the usual notation
            \begin{align*}
                \cos\sphericalangle(\myvec{x}, \myvec{y}) = \frac{\langle \myvec{x}, \myvec{y} \rangle_2}{\Vert \myvec{x} \Vert_2 \Vert \myvec{y} \Vert_2}, && \sin\sphericalangle(\myvec{x}, \myvec{y}) = \sqrt{1 - \cos^2\sphericalangle(\myvec{x}, \myvec{y})},
            \end{align*}
            where $\sphericalangle(\myvec{x}, \myvec{y})$ denotes the angle between two vectors $\myvec{x}, \myvec{y}\in \R^N.$

            Furthermore, we note that for both the $\ell_1$ and $\ell_2$ discrepancy terms, the centroids given in \cref{eq:l1:MinV:ArgMin} and \cref{eq:l2:MinV:ArgMin} cannot be computed directly from the data points $\{\myvec{x}_m \ \vert \ m=1,\dots,M \}$ and are dependent from the cluster membership matrix.

            Moreover, it is interesting to examine shortly some unregularized special cases and compare them with some results throughout the literature. Pompili et al.\ in \cite{Pompili:2014:CompMethod} consider the $\ell_2$ discrepancy term with $\lambda_{\myvec{U}} = \lambda_{\myvec{V}} = \mu_{\myvec{U}} = \mu_{\myvec{V}} = 0$ and further restricts $\myvec{U}$ to satisfy $\myvec{U}^\intercal\myvec{U}=\myvec{I}_{K\times K}$ by including the additional normalization constraint $\langle \myvec{u}_k,\myvec{u}_k \rangle_2 = 1 $ for all $k \in \{1,\dots,K\}.$ It can be shown, that the normalization constraint can be equivalently imposed on the rows of $\myvec{V}$ \cite{Pompili:2014:CompMethod}, which then leads to the distance measure $\dist(\myvec{x}_m, \myvec{v}_{\pi(m)})^2 = \Vert \myvec{x}_m - \langle \myvec{x}_m, \myvec{v}_{\pi(m)} \rangle_2 \myvec{v}_{\pi(m)} \Vert_2^2$ based on \cref{eq:l2:DistanceMeasure}. This coincides with the distance measure given in the proof of \cite{Pompili:2014:CompMethod}, where the equivalence to a weighted variant of the spherical $K$-means method is shown. For the centroids, we obtain $\myvec{v}_k = \nicefrac{1}{\Vert \myvec{u}_k \Vert_2^2} \sum_{m\in \mathcal{I}_k} U_{mk} \myvec{X}_{m,\bullet}$ based on \cref{eq:l2:MinV:ArgMin}.

            If we further restrict $\myvec{U}$ to be
            \begin{equation*}
                U_{mk} =        \begin{cases}
                                    0, & m\not\in \mathcal{I}_k, \\
                                    \frac{1}{\sqrt{\vert \mathcal{I}_k \vert}}, & m\in \mathcal{I}_k,
                                \end{cases}
            \end{equation*}
            which also leads to $\myvec{U}^\intercal\myvec{U}=\myvec{I}_{K\times K},$ we obtain the distance measure
            \begin{equation}\label{eq:l2:DistanceMeasure:SpecialCase:UOrthogonal}
                \dist(\myvec{x}_m, \myvec{v}_{\pi(m)})^2 = \Vert \myvec{x}_m - \frac{1}{\sqrt{\vert \mathcal{I}_{\pi(m)} \vert}} \myvec{v}_{\pi(m)} \Vert_2^2.
            \end{equation}
            by considering the minimization task \cref{eq:l2:U:MinTask} and inserting the non-zero solution $1/\sqrt{\vert \mathcal{I}_{\pi(m)} \vert}.$ Furthermore, based on \cref{eq:l2:MinV:ArgMin}, the centroids are given by
            \begin{equation}\label{eq:l2:Centroid:SpecialCase:UOrthogonal}
                \myvec{v}_k = \frac{1}{\sqrt{\vert \mathcal{I}_{k} \vert}} \sum_{m\in \mathcal{I}_k} \myvec{x}_m.
            \end{equation}
            Both the distance measure and the centroid given in \cref{eq:l2:DistanceMeasure:SpecialCase:UOrthogonal} and \cref{eq:l2:Centroid:SpecialCase:UOrthogonal} coincides with the ones stated in \cite{Pompili:2014:CompMethod}. If we instead constrain $\myvec{U}$ to be
            \begin{equation}\label{eq:Constraints:U:Classical K-means}
                \begin{aligned}
                    U_{mk} =        \begin{cases}
                        0, & m\not\in \mathcal{I}_k, \\
                        1, & m\in \mathcal{I}_k,
                    \end{cases}
                \end{aligned}
            \end{equation}
            we obtain based on the above procedure the usual Euclidean distance and the mean of the data points given by
            \begin{align*}
                \dist(\myvec{x}_m, \myvec{v}_{\pi(m)}) &= \Vert \myvec{x}_m - \myvec{v}_{\pi(m)} \Vert_2,\\
                \myvec{v}_k &= \frac{1}{\vert \mathcal{I}_{k} \vert} \sum_{m\in \mathcal{I}_k} \myvec{x}_m,
            \end{align*}
            used in the classical $K$-means approach. Finally, for the $\ell_1$ discrepancy term and the constraints in \cref{eq:Constraints:U:Classical K-means}, we obtain
            \begin{align*}
                \dist(\myvec{x}_m, \myvec{v}_{\pi(m)}) &= \Vert \myvec{x}_m - \myvec{v}_{\pi(m)} \Vert_1,\\
                V_{kn} &= \med (\tilde{\myvec{X}}_{\bullet,n}),
            \end{align*}
            based on \cref{eq:l1:DistanceMeasure} and \cref{eq:l1:MinV:ArgMin}, which yields the usual $K$-median approach with the $\ell_1$ distance and the median of the data points. Note that due to the application of the median, the centroids are chosen among the available data points $\{\myvec{x}_m \ \vert \ m=1,\dots,M \}.$ Furthermore, the application of the $\ell_1$ distance function leads to a more stable clustering approach with respect to outliers.
    \section{Conclusion} \label{sec:Conclusion}
        In this work, we presented a novel derivation framework to obtain connections between generalized $K$-means clustering approaches and regularized ONMF models by directly extracting the distance measures and centroids of the $K$-means method based on the considered ONMF problem. We applied this technique to non-standard ONMF models with elastic net regularization on both factorization matrices and derived the corresponding distance measures and centroids of the generalized $K$-means model. Furthermore, we gave an intuitive view on the obtained results and shortly described the effect of some regularization terms in the ONMF models on the distance function and the clustering outcome. Finally, we analyzed several special cases and found that the obtained results based on the proposed framework coincide with the results in the existing literature.

        Several further research directions could be of interest. A good starting point could be a further theoretical as well as numerical evaluation of the obtained generalized distance measures and centroids in \cref{sec:Regularized ONMF Model and Generalized K-means}. Moreover, a possible extension of this work could be to study whether the proposed derivation framework can be applied to more general discrepancy and regularization terms. Especially the analysis of gradient based regularization terms, which enfore spatial coherence in the clusterings, and their effects on the obtained distance measures and centroids of the generalized $K$-means approach could be of particular interest, since this can lead to an improved clustering performance for specific applications \cite{Fernsel:2021:ONMFTV}.
        
        A further theoretical direction constitutes the extension of the whole framework to infinite dimension spaces leading to continuous factorization problems. In combination with gradient based regularization terms, the analysis of first order conditions could lead to connections to $K$-means methods in a continuous setting as well as partial differential equations, whose solutions could lead to the desired centroids and distance measures.

    \section*{Acknowledgments}
    This project was funded by the Deutsche Forschungsgemeinschaft (DFG, German Research Foundation) within the framework of RTG ``$\pi^3$: Parameter Identification -- Analysis, Algorithms, Applications'' -- Project number 281474342/GRK2224/1.

    \bibliographystyle{siam}
    \bibliography{main}

\end{document}